\documentclass{article}
\usepackage[utf8]{inputenc}
\usepackage{amsmath,amssymb,amsfonts,latexsym}
\usepackage{graphicx}
\usepackage[all,cmtip]{xy}

\usepackage{times}
\usepackage{a4wide}

\usepackage{url}
\usepackage{hyperref}
\usepackage{xcolor}
\hypersetup{
    colorlinks,
    linkcolor={red!50!black},
    citecolor={blue!50!black},
    urlcolor={blue!80!black}
}
\usepackage{color,soul}

\usepackage[natbib=true,style=apa,backend=biber]{biblatex}

\title{Typed Mathematical Text for On-screen Examinations}
\author{Laura Kobel-Keller and Chris Sangwin}
\date{\today}

\begin{document}

\maketitle

\begin{abstract}
This paper discusses digital online mathematics examinations --- a discussion ranging from high school to university level examinations.  
In particular, we consider the nature of mathematical writing, what is distinctive about mathematical writing, and how mathematics can be typed into a machine.
This includes a review of features of notation and layout unique to mathematics and a survey of current technology for typed mathematics, including LaTeX and contemporary proof-checkers such as Lean.
Artificial intelligence has already been highly successful for optical character recognition, generating text from hand writing and is even increasingly applied to assess students' work itself.
A human-editable text-based format in the middle is important, but neglected.
The design of digital online mathematics examinations, which take this text as the source of truth for students' work, will have a profound effect on how mathematics is perceived and how mathematical activity is mediated.
Moving examinations on-screen effectively is an important design challenge and responsibility to future generations.
The challenge is to design software tools which support mathematics, that is tools which recede into the background and support the generation of mathematical work.
We argue for a human-editable text-based format, which includes semantic elements, at the heart of the process. 
\end{abstract}


\section{Introduction}

This paper discusses digital on-screen mathematics examinations.
Finland has already fully digitized its high school graduation examinations, including mathematics, which have been conducted online via the Abitti platform since the spring of 2019. 
In Finland, students use their own laptops to connect to a local server in school, covering all subjects in a secure, digital format, replacing paper tests with modern digital assessments.
The International Baccalaureate Organisation have started a major project to move all their exams, including mathematics, on-screen by 2030.
Section \ref{sec:eth} gives details of the examination system in use at ETH Z\"urich currently.
At the start of 2026 it is clear that more mathematics examinations will soon take place on-screen, both in school and undergraduate examinations.

Controlled conditions provided by invigilated examinations (both written and oral/viva) is one, and perhaps the only, reliable way to distinguish students' work from the output of large language models or other artificial intelligence (AI) tools and to assess \emph{individual} performance.
Impersonation of one person by another has always been a form of cheating in examinations.
The availability of AI provides a readily accessible opportunity for cheating by all students, particularly with unsupervised coursework.
This paper does not discuss cheating in examinations (in any form), alternative forms of assessment, or similar issues, see \citep{Ellis2024}.
Instead, we discuss the nature of mathematical activity generally.
How written examinations can support students' efforts to learn mathematics, and indeed the case {\em for} examinations, has been made elsewhere, e.g.~\citet{Brereton1944}.
This paper assumes that students will continue to learn to write correct mathematics themselves.
This paper also assumes that students will write to learn, as they have always done, practicing classical techniques and seeking to find meaning and understanding alongside fluency and mastery of the learnt content.
This paper discusses options, available in 2026, to support the practical move to online mathematics examinations and exercises.
We end with a brief discussion of how AI technology might change what students learn in Section \ref{sec:future}.

For many essay subjects the move online is entirely natural.
Indeed, the only time most humanities students write with a pen or pencil for extended periods of time is during an examination.
Mathematics is different.
The continuing popularity of pen, paper and chalk for mathematical work is not conservatism \citep{Greiffenhagen2014}.
Writing really is an important part of professional mathematical practice and thinking.
Mathematical notation, diagrams and pictures are all essential to mathematical cognition.
The challenge is to design software tools which support mathematics: a genuine tool recedes into the background of activity.
Generally speaking, tools should enhance, not distort, cognitive processes involved in any activity.
Indeed, work such as \textcite{Zhang1994} suggest that the external environment mediates cognition, and that information is distributed both internally and externally.
\begin{enumerate}
  \item Tools must mediate activity effectively; poorly designed tools distort or obstruct activity or even tempt users to be lazy and outsource the actual process of thinking.
  \item Tools used for high-stakes assessments must be identical to those used during learning.
\end{enumerate}

One option for supporting on-screen examinations is to scan and digitise writing on paper or with a graphics tablet.
Optical character recognition (OCR) for mathematics is now excellent.
Reliable tools, such as Mathpix (see {\tt https://mathpix.com/}, Mar 2026), take an image and produce machine readable, searchable text (e.g.~in LaTeX format).
These tools are very useful and could be used in examinations.
This article does {\em not} discuss the use of OCR in examinations, rather this article looks at the nature of mathematical text and concentrates on typed mathematical text.
We include observations, in the context of AI supported grading, that OCR processing might not be necessary.

\section{Digital examinations at ETH Z\"urich}\label{sec:eth}

The Federal Institute of Technology Zurich (ETH Z\"urich), is a world leading research and teaching institution and one of the largest public universities in Switzerland that focuses on (natural) science, engineering, technology and mathematics (STEM), with over twenty five thousand students in 2025.

All students at ETH have to take several courses in mathematics and so the Department of Mathematics has one of the highest load of examinations.
Almost all the mathematics examinations for students in study programs other than mathematics are written, and so are the first examinations for students in mathematics itself. 
Oral exams appear only during the second study year for students in mathematics. 
As a consequence, several hundred written exams are taken in each session.
In the last 20 years, the number of students has doubled, and another doubling is expected. 
This second doubling is due to an expected general increase in student numbers and the planned change of the academic calender.
There are also some major structural changes, one of which is the change from year-round courses to exclusively courses with examinations at the end of each semester. 
Another structural change comprises the shift from one single, final examination, to a combination of several smaller assessments during the semester with a smaller final examination.


As one of the leading educational institutions worldwide, ETH strives to implement the highest quality assessments, including high quality grading processes.
The need for assessment integrity means that invigilated assessments, much like traditional examinations, will remain an important component for the foreseeable future. 
The aim of moving to digital examinations --- or at least to establish such exams as one of the standard formats --- is to maintain the high quality of traditional examinations, and to develop a new level setting international standards.
Rather than replicate paper-based examinations in a digital format, the goal is to provide modern digital examinations which include giving students access to materials, like text, data material, video clips, and specialist software such as code development environments seeded by a non-trivial starting point.
The ultimate goal is to incorporate more examination exercises that go beyond the classical questions where standard recipes have to be applied, with problems close to real applications of the sort that the students will encounter later on in their professional life. 
From a purely didactic point of view, such problem solving during an examination should be enhanced --- albeit accepting all the challenges that come along with such a change. 

Traditionally, all written exams in mathematics contained several classical open questions (e.g. optimization problems, applications of integral theorems (Gauss, Green, Stokes)).
A first change in the ETH examamination paradigm was the incorporation of multiple choice questions. 
This started more than 10 years ago for examination in non-mathematics study programs. 
For the grading of multiple choice questions at ETH, we use the automated grading system AMC (see \url{https://www.auto-multiple-choice.net}) with several modification and adaptions to our specific needs. The AMC system allows for shuffling potential answers and questions and is also one of the technical pillars for exams et EPFL Lausanne.
Although controversial from a didactic point of view, in the last year the established standard is to cover about half of the total possible points in an examination by such multiple choice questions or box questions where only the final answer is considered (e.g. the result of a moderately difficult integral computation) --- possibly with the option to give partial credit.  

For examinations in the study program in mathematics and physics, multiple choice questions were also introduced but only as recommended and at a lower percentage. The bigger change in paradigm for those examinations was the introduction of exercises in which the students could reproduce mathematical proofs that they saw during lectures and that the students were required to learn by heart.

Our first experience with digital exams in 2020 was with non-compulsory mid-term tests, mainly containing multiple choice questions, and shortly after with final compulsory exams. 
Concerning the bigger compulsory exams, the first designs contained multiple-choice questions, box questions where only the final answer was graded and typically also one or two open questions where the questions were displayed on-screen but where the answers had to be written classically on paper.

Now we have significant diversity in digital exams.  The simplest digital exams consist of only only multiple choice questions which are solved on-screen. 
In more complex exams, the students have to code, or to show other digital skills, and for these we use our in-house-developed platform ``code expert'' (see \url{https://learning-teaching-fair-2024.ethz.ch/project/codeexpert/}) with Jupyter notebooks. This diversity also reflects the different needs, didactic principles, and preferences of lectures.  

One unifying aspect of all variants of digital exams is that ETH uses the learning managing platform Moodle in combination with the in-house-developed Safe Exam Browser (see \url{https://safeexambrowser.org/download_en.html}).
Our Moodle based system incorporates STACK questions \cite{2013CAA}. 
Students use bring-your-own-device as well as managed devices from ETH, the first for non-compulsory assessments the latter for compulsory exams.

In order to further bridge all the various settings and requirements, gain in efficiency and on-site data control, ETH identified the desire for one unified online platform where all types of assessments, both during the semester as well as final exams, and all types of questions (open, multiple choice, short answer, drag and drop questions, coding etc.) and all grading options.

Traditionally, mathematics examinations assess both calculation and reasoning skills, see \citep{Smith1996}. 
In order to incorporate both competencies, the current standard practice is a combination of multiple choice questions; closed questions where only the final result is graded; and open questions in which students solve more elaborate problems, justify each step or give proofs of mathematical statements. 
There are different technical solutions for digital grading including AMC (for multiple choice questions), STACK, and other automatically graded Moodle question types.
Automatic grading enables fast, efficient, grading without any bias, and provide excellent options for questions that have either a numerical result, a mathematical expression or a single keyword as solution or options where scaffolded exercises are desired.
That said, creating robust, high quality, autocorrected questions is still challenging technical work.

The capability to put logical pieces of information correctly together, e.g.~giving a mathematical proof, can still only partially be automatically assessed, see \citep{2023-Proof-Assessment}.
Thus, from a didactic point of view, even with the best solution at hand for automatically corrected questions, ETH sees no end for a need for open questions.
Unfortunately, no good solution exists for open questions, which still have to be graded manually, even when using very efficient grading platforms.
Unlike other disciplines where typed text is sufficient, open questions in mathematics still have to be written on paper, scanned and uploaded. 
So, the current situation where whatever the students write on paper has to be digitalized first is unsatisfactory. 
Even in the context of AI assisted grading, an intermediate step where the students could enter their mathematical expressions in a digital manner instead of handwriting would be most beneficial, both for reliability and efficiency.
A human-readable digital format could be provided to students for editing, or typed directly (or indeed both).
Ideally, what is needed is a way for students to type their working and submit a complete mathematical argument with the final answer as a digital document.
Hence, the focus of this paper is a format for typing complete mathematical arguments.
The next sections discuss the nature of mathematical writing before we discuss options for such a typed format for mathematics.

\section{The nature of mathematical writing}

This section provides some examples of features specific to mathematical writing and discusses the nature of mathematical work.
Basic notational features needed for written mathematics include the following.
\begin{itemize}
  \item Numerous special symbols, e.g. $\times,\cup,\in,\infty$, which do not appear on keyboards.
  \item Regular use of superscripts (e.g.~$x^2$) and subscripts (e.g.~$a_1$).
  \item Use of a variety of character sets, especially Greek (e.g.~$\lambda$) and ``Blackboard bold''  (e.g.~$\mathbb{R}$).
  \item Use of a vinculum (a long line) to group terms, e.g.~in fractions \(\frac{1}{1+x^2}\) and surds \(\sqrt{a^2+b^2}\).
  \item Matrices, tables and other two-dimensional layout.
\end{itemize}

Mathematics often contains diagrams, such as the commutativity diagram below, which shows how the operations of addition \(a+b\) and multiplication \(a \times b\) are related.
\begin{equation*}
\xymatrixcolsep{3pc}
\xymatrix{
   a        \ar@<2pt>[r]^{e^x} \ar[d]_{+} & A          \ar@<2pt>[l]^{\ln(x)} \ar[d]^{\times} \\
   a+b      \ar@<2pt>[r]^{e^x}            & A\times B  \ar@<2pt>[l]^{\ln(x)} 
}
\end{equation*}

Elementary Mathematics contains two important aspects: geometry and computation.
Descartes linked the two by using coordinates and applying algebra.
Modern mathematics emphasises the concept of functions, and an important topic in elementary mathematics is the study of graphs and calculus.
\citet{Brown2008} argues persuasively that rigour can be found in some diagrams: ``{\em I don’t see any abandoning of rigour by allowing the legitimacy of picture-proofs.}'' \citep[p.~46]{Brown2008}.
However, typically in education students are expected to provide written, formalized, proofs.

Mathematics often makes use of layout to structure complete mathematical arguments.
The most important example in school mathematics is {\em reasoning by equivalence} \citep{2018Sangwin-equivalence-proof}.
For example, a quadratic equation can be solved as follows.
\[
\begin{array}{llll} 
                 & x^2-10\,x+9=0                            &  \cr 
\Leftrightarrow  & {\left(x-5\right)}^2-16=0                & \mbox{Complete the square} \cr 
\Leftrightarrow  & \left(x-5+4\right)\left(x-5-4\right)=0   & \mbox{Difference of two squares noting } 16=4^2\cr 
\Leftrightarrow  & \left(x-1\right)\left(x-9\right)=0       & \cr  
\Leftrightarrow  & x=1\,{\text{ or }}\, x=9                 & (AB=0 \Leftrightarrow A=0 \mbox{ or } B=0)\cr 
\end{array}
\]
Reasoning by equivalence extends naturally to work with inequalities and calculus where equivalence is not retained between lines.
Many of these arguments include both detailed steps in the calculation and justification of the legitimacy of those steps.
An essential detail in mathematical work is the ability to {\em cross reference} between steps in an argument or between different parts of an argument, and any nested sub-arguments.
In examinations, justification is typically essential and marks for accuracy are often only awarded when appropriate justification is given by the student.  \citep{2016AutomationExaminations} provides some examples and data.

\begin{figure}[ht]
\begin{center}
\includegraphics[width=7cm]{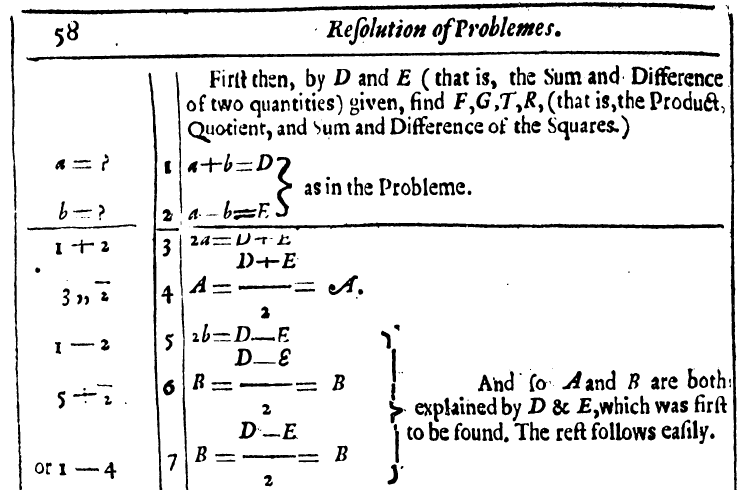}
\end{center}
\caption{An algebraic derivation in two-column format}\label{fig:Pell1668}
\end{figure}

\begin{figure}
\begin{center}
\hfill
\includegraphics[width=7cm]{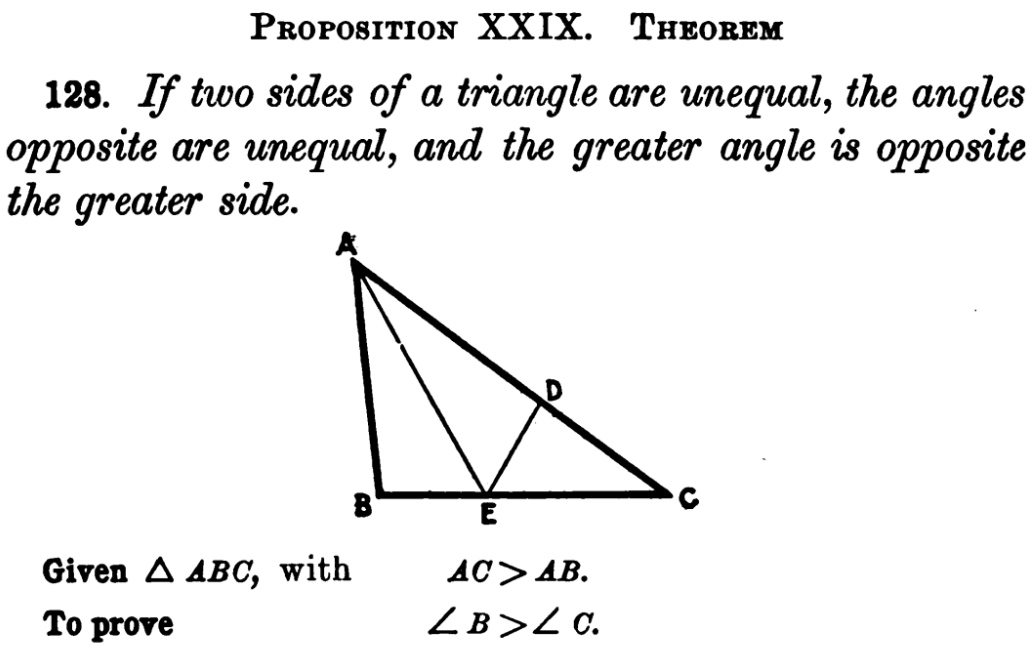}
\hfill
\includegraphics[width=7cm]{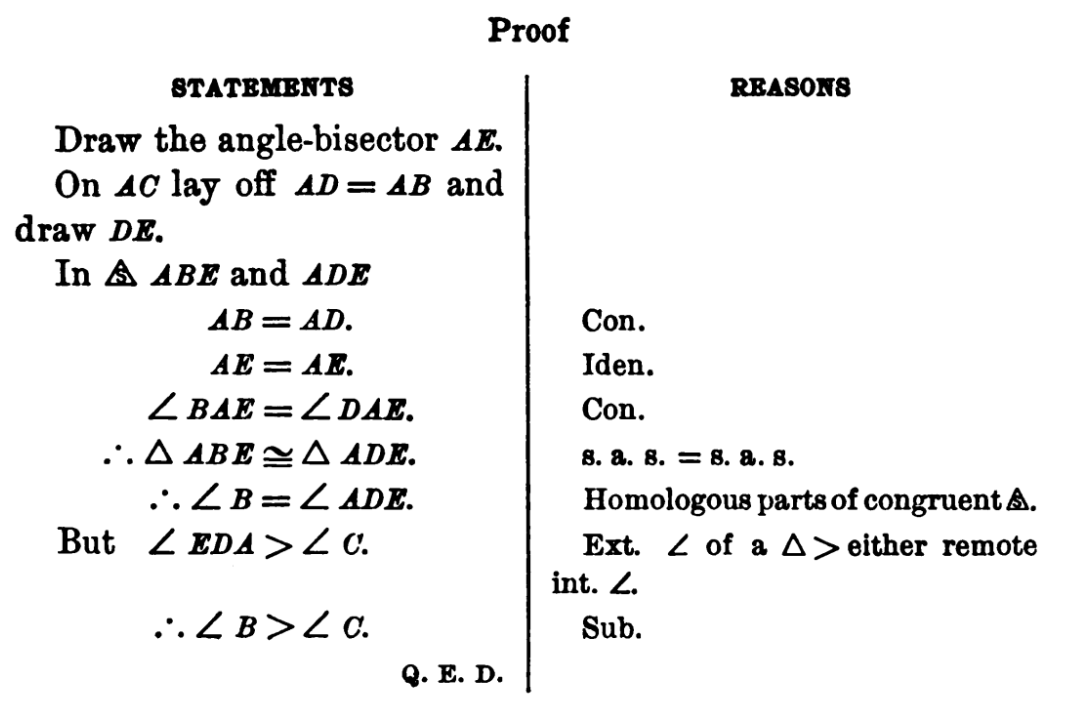}
\hfill
\end{center}
\caption{A school geometry proof in two-column format}\label{fig:Schwamb1904-pg67}
\end{figure}

Two-dimensional layout can be formalised, and laying out arguments in two-columns can be traced back hundreds of years.
An algebraic example from  \citet{Brancker1668} is shown in Figure \ref{fig:Pell1668} (which contains somewhat archaic algebraic $17^{th}$ Century notation).
The two-column format is enduring popular.
An example of a two-column proof, from  \citet{Schulze1913}, p.~67, is shown in Figure \ref{fig:Schwamb1904-pg67} with the theorem to the left, and the two-column proof to the right.
As with all things, when a particular format becomes routine, meaning can be lost and so the two-column proof format has been criticised, e.g.~see Weiss et al. \citep{Weiss2009}.
A recent discussion of mathematical style was given by \citet{Ording2019}, listing numerous ways of laying out mathematical arguments.
Providing a structured interface on-screen is possible, and has some advantages, and we provide an example in Section \ref{sec:Back}.

\begin{figure}
\begin{center}
\includegraphics[width=8cm]{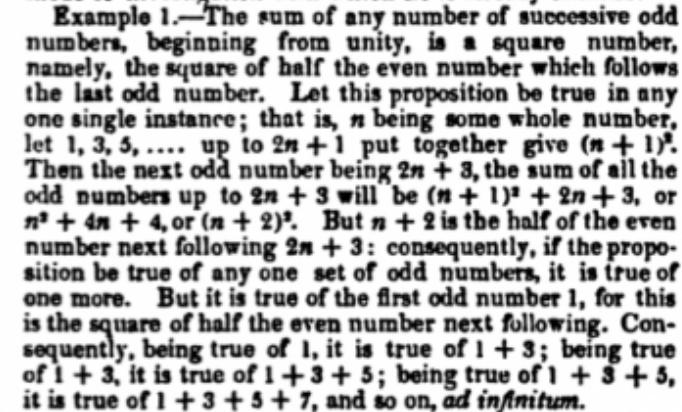}	
\end{center}
\caption{De Morgan's proof that \(\sum_{k=1}^n (2k-1)=n^2\)}\label{fig:DeMorgan1836-induction}
\end{figure}

Mathematical arguments {\em could} be written with less notation, and in plain paragraph style.
An early, if not the first, use of induction is found in \citet{DeMorgan1838}.
His argument that the sum of the first $n$ odd numbers is $n^2$ is given in Figure \ref{fig:DeMorgan1836-induction}.
Notice how difficult this dense text is to read (regardless of the reproduction quality or typeface size), compared to a modern induction proof which uses more algebraic calculation, and is laid out with more space on the page.
See Figure \ref{fig:2025-unicode} for a more contemporary version, and \citep{2023-Sum-Odd-Num} for a discussion.
On the other hand, notation can become excessive: effective and judicious use of notation is something of a skill in itself.
Mathematicians have long acknowledged this.
\begin{quote}
``The advantage of selecting in our signs, those which have some resemblance to, or which from some circumstance are associated in
the mind with the thing signified, has scarcely been stated with sufficient force: the fatigue, from which such an arrangement saves
the reader, is very advantageous to the more complete devotion of his attention to the subject examined.", Babbage
 \citep[p.~370]{Babbage1827}
\end{quote}
The historic view is reinforced by modern theory of cognitive science, e.g. ~\textcite{Zhang1994}: the external environment mediates cognition.
Notation really does matter.

\subsection*{History of, and innovations in, notation}

Shared notation and shared meaning is essential, without which effective communication is impossible.
For example, the complex logarithm is a multivalued inverse of \(w\rightarrow e^w\).
The multivalued inverse of the function \( w \rightarrow we^w\) is known as the Lambert $W$ function. 
\[ W(x)=y\quad \Leftrightarrow\quad x=ye^y.\]
The $W$ function was first used around 1758 and  \citet{Corless1996} conclude as follows. {\em ``Names are important. The Lambert W function has been widely used in many fields, but because of differing notation and the absence of a standard name, awareness of the function was not as high as it should have been."}

Evolution of notation is very slow, and mathematics is conservative, see e.g. Cajori \citep{Cajori1928}.
Furthermore, it appears very difficult for an individual to change the notation they are used to using.
As an example, the logarithm of $x$ to base $a$ is often written \(\log_a(x)\).
\citet{Brown1974} suggested that the notation \(a\downarrow b\) be used instead.  
Rather than $\downarrow$, we use subscripts to write \(a_b\):  
\[a^x=b \quad \Leftrightarrow \quad x = a_b.\]
Then, \( a_{b\times c} = a_b + a_c\), together with the rules:
\[ a_{b^c} = (a_b)\times c,
 \quad (a_b)\times (b_a) = 1,
 \quad (a_b) = (a_c)\times (c_b).
\]
These symbolic rules have what is called {\em visual salience}.
\begin{quote}
    ``The quality of visual salience is easy to recognize but difficult to define.  Visually salient rules have a visual coherence that makes the left- and right-hand sides of the equations appear naturally related to each other. ", see Kirshner and Awtry \autocite{Kirshner2004}
\end{quote}
Aesthetics play an important role in mathematics, and visually salient rules have an important aesthetic appeal.
Furthermore, using identical or similar notations enables effective communication via analogy, cf. Polya \citep{Polya54}.
An important example is the use of exponential notation $x^n$ to denote the repetition of an operation, such as multiplication, $n$ times.  The following theorem is then an immediate consequence of the notation, albeit with an underlying induction.
\begin{equation}\label{th:powers}
    \mbox{If }n,m\in\mathbb{N} \mbox{ then } a^na^m=a^{n+m}.
\end{equation}
Exponential notation is re-used regularly, e.g.~to denote repeated differentiation.
Such coherence, elegance of expression, and effectiveness is very appealing and some mathematicians were explicit in their desire to retain and exploit this coherence.
Indeed, it became a conscious meta-theorem: the {\em principle of permanence} (sometimes the {\em law of the permanence of equivalent forms}), promoted by mathematicians such as George Peakcock \citep{Peacock1842v1,Peacock1842v2}, see \citet{Lambert2013}.

The principle of permanence demands that the theorem in \eqref{th:powers} should remain and become a {\em definition} when the domain for \(n,m\) is extended, to the rational numbers say.
This permanence {\em requires} \(a^{\frac{1}{2}}\) to be the square root.
Furthermore, this principle of permanence was highly influential in understanding fractional calculus, e.g.~the meaning of \(\frac{\mathrm{d}^{n}}{\mathrm{d}x^{n}}\), when \(n=\frac{1}{2}\) see \citet{Debnath2004}.

Another example of effective notation choice is writing $x^{\underline{n}}$ as the $n$th falling product of $x$.
Just as \(x^n\) is notation for {\em repeated multiplication} define
\(x^{\underline{n}}\) as notation for the {\em falling product}
\begin{align*}
  x^{\underline{2}} & = x(x-1) \\
  x^{\underline{3}} & = x(x-1)(x-2)\\
  \vdots           & \\
  x^{\underline{n}} & = \underbrace{x(x-1)(x-2) \cdots (x-(n-1))}_{n \mbox{ terms}}
\end{align*}

There is a profound analogy between calculus of a continuous single-variable function, and summation of finite discrete series.
This analogy can be expressed very clearly with a judicious choice of notation, such as the falling product.
For example
\begin{equation}\label{eq:sum_falling_prod}
  \sum_{k=0}^{n-1} k^{\underline{m}} = \frac{1}{m+1} n^{\underline{m+1}}.
\end{equation}
In the case $m=3$ we have
\[
  \sum_{k=0}^{n-1} k(k-1)(k-2) = \frac{1}{4} n(n-1)(n-2)(n-3)
\mbox{ or, written another way, }
  \sum_{k=0}^{n-1} k^{\underline{3}} = \frac{1}{4} n^{\underline{4}}.
\]
With this notation one can't help notice the similarity with the continuous counterpart:
\[ \int x^m\mathrm{d}x = \frac{1}{m+1} x^{m+1},\quad (m\neq -1).\]
Appreciation of the systematic correspondence between the traditional calculus and the method of finite differences is seriously hampered by a lack of agreed notation and terminology, e.g.~that for falling products.
\citet[p.~81]{Boole2007} complained about this in the 1870s.
\begin{quote}
``The student will doubtless already have perceived how much the branch of mathematics that forms the subject of our present consideration suffers from it's not possessing a clear and independent set of technical terms. [...] There is no reason why the present state of confusion should be permanent, [...]" 
\citet[p.~81]{Boole2007}
\end{quote}

All the basic topics (summing finite series, calculus, etc.) are already included individually, but separately, in school mathematics however the analogy is rarely explored explicitly or systematically.
Lack of suitable notation gets in the way.

This historical interlude serves to reinforce the importance of notation, both to the subject as a whole and the role notation plays in the cognition of an individual mathematician.
Indeed some authors, such as \citet{Babbage1830}, have treated mathematical notation as an explicit design problem, as did \citet{Dijkstra2002}.
However, for better or for worse, online examinations need to support traditional algebraic notation in its current form.

\section{Mathematics and computer typesetting}

This section records, for completeness, details of the most important contemporary techniques for typing and encoding mathematical text using English, or a left-to-right European language, using a computer.

\subsection{Native ASCII and unicode text}

The American Standard Code for Information Interchange (ASCII) is an encoding of 95 printable and 33 control characters (128 in total).
For example, capital letter $A$ is code 65 and the space is code 32.
The ASCII standard was first published in 1963, was revised and extended numerous times, and remains influential to this day.
ASCII is the basis for most typed text, indeed most contemporary keyboards have about 60 primary keys, with perhaps about 100 keys including function keys and numeric keypads.
This limits the number of symbols available to be typed directly.
While some other symbols are available, e.g.~the Euro symbol is commonly available on UK keyboards using {\tt AltGr+4} in Windows, only a few of the commonly used mathematical symbols can be typed directly.

\begin{figure}
\begin{center}
\includegraphics[width=10cm]{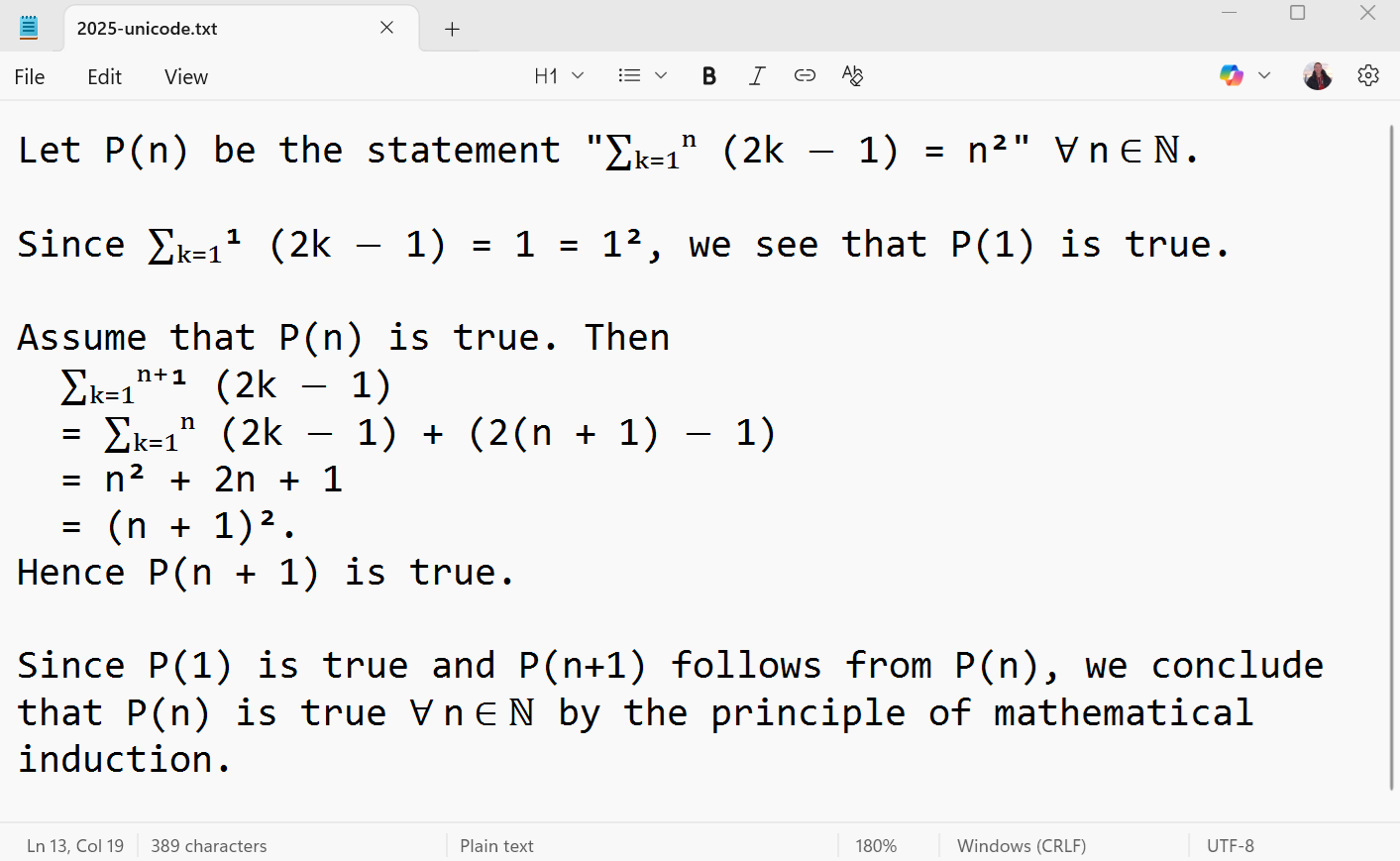}
\end{center}
\caption{Editing a unicode document using the Windows ``Notepad'' text-editor}\label{fig:2025-unicode}
\end{figure}

Contemporary computers provide much richer native support for symbols and fonts as unicode text, regardless of whether they can be typed directly.
The Unicode Standard includes a significant proportion of symbols used in mathematics, and virtually all those needed for school mathematics.
Unicode contains dedicated symbols for superscript and subscripts.
Figure \ref{fig:2025-unicode} shows an example of native unicode text.
Just to be clear, this is a plain text file ({\tt .txt} file) opened in Windows Notepad.
Unicode has no problems typesetting \(\forall n\in \mathbb{N}\) but the superscripts and subscripts in \(\sum_{k=1}^n\) are not aligned perfectly (i.e.~\({\sum_{k=1}}^n\)) because there is no way to stack symbols vertically in a traditional way.
Unicode is useful for very simple mathematics, especially simple expressions with single-level superscripts or without complex two dimensional layout.
However, typing unicode mathematics characters is difficult.  There are a number of options.
\begin{enumerate}
  \item Select Unicode characters visually from a palette, such as the Character Map on Windows.
  \item Type in the code explicitly.  For example, on some computers holding down the {\tt Alt} key and typing a code (e.g.~0247) yields the corresponding unicode symbol (e.g.~$\div$).
  \item Cut and paste from elsewhere, convert one format to another, (or use OCR).
\end{enumerate}
None of these options are satisfactory for examinations as they interrupt the train of thought and are a distraction.
Nevertheless, unicode text is supported almost universally and is very simple to implement cross platform.

\subsection{AsciiMath and Space Math}

AsciiMath is a markup language for mathematics \citep{Cervone2012}.
The goal is that users type only ASCII characters, or characters available on all keyboards, to create complex mathematical expressions.
Rather than have a dedicated key for the sigma/sum character $\sum$, users type the text-string \verb$sum$.
For example, a user types 
\begin{verbatim}
sum_(i=1)^n i^3=((n(n+1))/2)^2
\end{verbatim}
which is displayed as 
\[ \sum_{i=1}^n i^3=\left(\frac{n(n+1)}{2}\right)^2\]
See {\tt https://asciimath.org/} (May 2026) for technical details.
There are many advantages to AsciiMath.
\begin{itemize}
  \item Simple, direct keyboard entry, without too many distractions.
  \item Markup codes have names corresponding to their (English) operations, e.g.~\verb$sum$ or \verb$beta$ (Greek letter name).  
   These codes can be learned and follow a somewhat predictable pattern.
  \item Use of AsciiMath results in a clean text-document which is completely transparent, and which can include unicode characters as well.
  \item Matrix entry is relatively unobtrusive, e.g.~\verb$[[a,b],[c,d]]$.
  \item In larger documents, users typically enclose mathematics using back-ticks \verb$`...`$ to distinguish mathematics from text. 
  \begin{verbatim}
  This is text with an equation: `a^2=b^2+c^2`.
  \end{verbatim}
\end{itemize}
However, there are also non negligible disadvantages of AsciiMath, include the following.
\begin{itemize}
  \item Not WYSIWYG: a separate interpreted display panel is needed.
  \item Some learning of codes is needed.
  \item A limited and fixed set of mathematical operations is supported.
  \item Limited attempt to encode semantic meaning, preserving well-known ambiguities in mathematical notation.
\end{itemize}


Space Math builds on ideas from AsciiMath to address the problem that current ways of writing mathematics do not provide an easy way to distinguish between common ambiguous expressions.  
For example, in \(x(t+1)\) it is not clear if the author intended $x$ to be a function  evaluated at the point $t+1$, or a variable multiplied by $t+1$.  The original author had a context, but typing an expression such as \verb$x(t+1)$ is ambiguous.  Space Math allows semantics by the use of spaces, e.g. \verb$x(t+1)$ is interpreted as function application whereas \verb$x (t+1)$ is interpreted as multiplication.  Decisions can be made later how this difference is displayed (if at all).

Space Math also takes inspiration from the Python programming language to use whitespace for multi-line systems of equations for more advanced users.
\begin{quote}
  ``Space Math syntax incorporates ingredients from LaTeX, AsciiMath, and UnicodeMath. To the extent possible, Space Math is backward compatible with all three of those languages, meaning that the majority of expressions in those systems are parsed correctly when incorporated into Space Math".\\
  \verb#https://aimath.org/~farmer/spacemath/# (Jan 2026)
\end{quote}

If, during an exam, students forget to use the semantic aspects correctly and (for example) mix \verb$x(t+1)$ and \verb$x (t+1)$ it is very unlikely to cause any serious confusion.  
The examination provides a context, and examiners (for whom correct use of the semantic markup is not a learning outcome) can easily condone this ambiguity and understand the work {\em as they routinely do currently with handwriting}.

\subsection{TeX and LaTeX} \label{sec:tex}

TeX and LaTeX are the most popular typesetting systems for mathematics.
TeX was developed by Donald Knuth, with the first version available in 1978 \citep{KnuthTeX}.
Knuth's original motivation was to replicate traditional movable-type printing in electronic form, as an aid to physical print publication.
He was trying to typeset large structured mathematical documents, such as books.

Commands in TeX, including mathematics expressions, start with a backslash: $\backslash$.  For example, a user types 
\begin{verbatim}
\sum_{i=1}^n i^3=\left(\frac{n(n+1)}{2}\right)^2
\end{verbatim}
which is displayed as 
\[ 
\sum_{i=1}^n i^3=\left(\frac{n(n+1)}{2}\right)^2
\]
Notice here the command \verb$\left(...\right)$ are typesetting commands to automatically re-size the parentheses, and curly braces are used to group terms together.

TeX subsequently developed into LaTeX with a very wide range of optional packages for bibliography management, typesetting diagrams and supporting just about every other typesetting requirement.
TeX and LaTeX are difficult to learn.
Knuth explicitly, and by design, separated out the processes of (i) authoring (perhaps writing on paper); (ii) typesetting in TeX; (iii) processing the document automatically (e.g.~to create cross referencing) and finally (iv) printing (or viewing a .pdf on-screen).
By design TeX and LaTeX are not WYSIWYG.
The difficulty of learning TeX, and the nature of the commands needed, mean most people find it's quite distracting to typeset a complete argument in TeX as they derive it.
The result is beautiful, and still represents the best available system for typesetting mathematics, but also make it unsuitable for examination answers.
For many, TeX markup is too complicated and too distracting for deriving short, focused mathematical arguments without an immediate preview.

Furthermore, there is no suggestion that when developing TeX that Knuth wanted to encode the mathematical meaning of a text, so that further computation or proof verification could take place later.
Indeed, interpreting meaning from TeX source is well-known to be very difficult, e.g.~\citet{Fateman1999}.
In an examination such verification would include automatic assessment.
This is one example of where new technology {\em looks backwards}, replicating existing practices and systems.
The same happened with Babbage's mechanical computers.
Babbage's motivation was to accurately calculate and print tables of logarithms, replacing human-computed tables which were error-prone \citep[p.~138]{Babbage1864}.
We can find nothing which suggests Babbage anticipated everyone having their own computer, to help them compute directly.
It is entirely natural that the first use of technology is initially within traditional modes of working, teaching and for examinations.

Originally, TeX documents were compiled to eventually produce a printed document.
There are now many ways to use TeX/LaTeX more interactively, with an immediate preview.
There are also many online WYSIWYG editors that use TeX/LaTeX for mathematics, e.g.~{\tt https://www.lyx.org/} (Feb 2026) is particularly popular.

\subsection{Lean}

Lean is a proof assistant.
A proof assistant is software that checks proofs are formally correct.
Lean is a programming language enabling users to define mathematical objects, specify properties of these objects, and encode mathematical statements about these objects.
Proof assistants, by design, encode the mathematical meaning, not just typeset traditional mathematical documents for printing.

%
%

\begin{figure}
\begin{center}
\includegraphics[width=12cm]{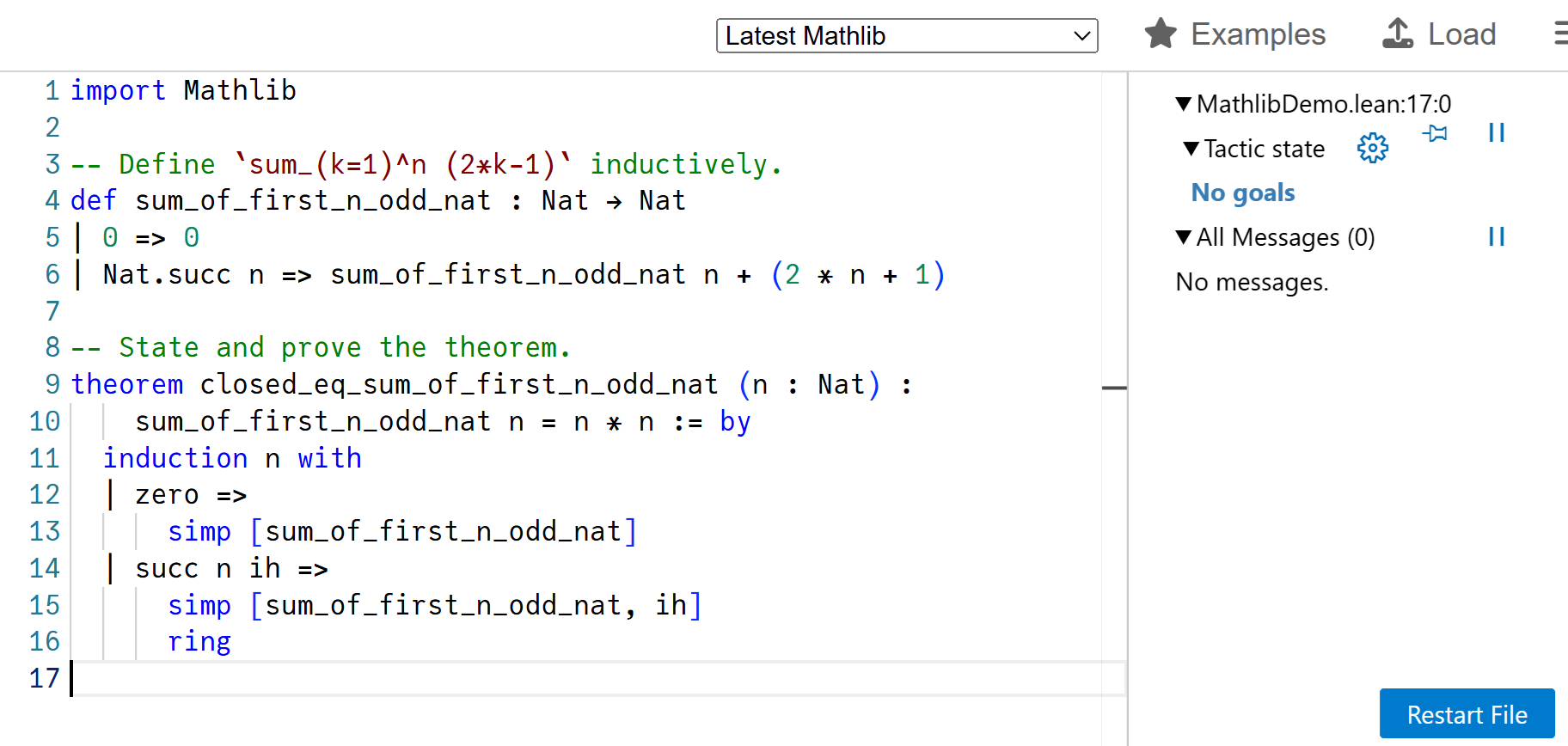}	
\end{center}
\caption{A proof verified using LEAN}\label{fig:Lean2025-induction}
\end{figure}

An example using the web-interface to Lean 4 at {\tt https://live.lean-lang.org/} (May 2026) is shown in Figure \ref{fig:Lean2025-induction}.
This example shows one way to verify that the sum of the first $n$ odd numbers equals $n^2$, the same theorem proved by De Morgan in Figure \ref{fig:DeMorgan1836-induction}.
Lines 4--6 inductively define the formula $\sum_{k=1}^n (2k-1)$ and lines 9--16 prove the theorem.
The theorem has been established on the right because the system reports ``no goals''.
The proof is complete when the set of outstanding goals is empty.

This example is included as it illustrates some very useful interface features.
\begin{enumerate}
  \item The student types a unicode text document, which can include special symbols.  E.g.~on line 4 of Figure \ref{fig:Lean2025-induction} contains the arrow symbol, $\rightarrow$, which does not appear on a keyboard.
  \item The interface is split, with a text-area for typed input on the left, and immediate feedback on the right.
  \item Some parts of the text have actual meaning, and comments can be added, e.g.~lines 3 and 8.
\end{enumerate}
A very similar two-panel interface is provided by {\tt https://www.overleaf.com/} (May 2026) for editing shared LaTeX documents online.

Adoption of proof assistants by the mathematical community is still gradual, and the use of proof assistants is by no means universally accepted, see e.g. MacKenzie \citep{MacKenzie2004}.
However, they are gaining users and there are notable successes in research, and increasingly in education. cf. Wemmenhove \citep{Wemmenhove2025}.
For every-day mathematics Lean is problematic.
Lean is a complex, and sophisticated, programming language which is difficult to learn.
For example, users need to remember commands such as {\tt ring}, which is used on line 16.
More seriously, by design, proof assistants hide a lot of details.
For example, compare line 16 in Figure \ref{fig:Lean2025-induction} with the four lines of algebraic work in Figure \ref{fig:2025-unicode} which prove $P(n)\Rightarrow P(n+1)$.
The use of the command {\tt ring} in Lean hides the traditional algebraic work, much of which is precisely what is tested by current examinations.
The Lean community are well aware of these issues, and project such as {\em Lean verbose} and {\em Waterproof} by Wemmenhove \citep{Wemmenhove2025} seek to address them.


\subsection{Equation editors}

\begin{figure}
  \begin{center}
  \includegraphics[width=0.6\textwidth]{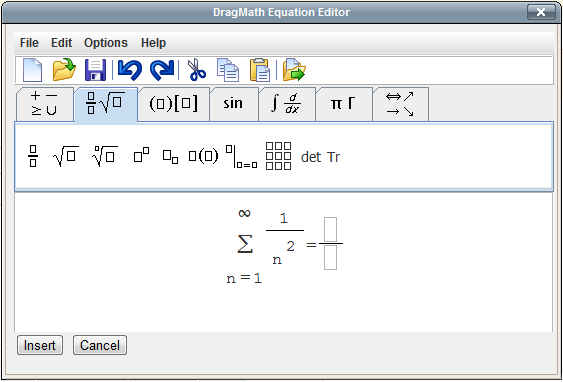}
  \end{center}
  \caption{The DragMath equation editor}\label{fig:DragMath1}
\end{figure}

Drag and drop equation editors are enduringly popular.
An example, DragMath, is shown in Figure \ref{fig:DragMath1} as described in \citet{2012DragMath}.
Editors help users build up a complex expression without having to remember commands or codes.
Users immediately see the expression as they build it up.
More experienced users typically find that switching from typing at a keyboard to clicking/dragging a mouse interrupts the train of thought, as the focus of attention switches back and forth between a text document and the equation editor.
This really is a serious problem for experienced users: learning to type certainly takes time.
Then again, learning to write with a pen or to use algebra effectively also take significant time and effort.
The goal of a frictionless typing system for mathematics, with preview, acknowledges the need to expend effort in learning and focused practice.

Every equation editor is slightly different, and if proprietary editors are used in examinations then they need to be readily available for learning.
Equation editors have proved very short-lived, and often platform dependent.
For example, older Word documents which used earlier MS equation editors are effectively un-readable now, and DragMath stopped working after Java applets required security certificates.
Equation editors often lack accessibility features, or cannot take advantage of accessibility tools (such as screen-readers).
Text-based formats do not suffer from these problems.

\begin{figure}
  \begin{center}
  \includegraphics[width=0.7\textwidth]{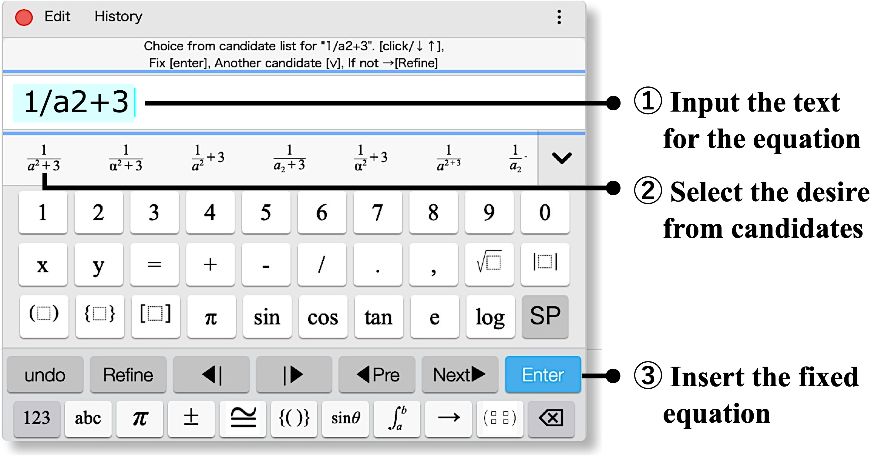}
  \end{center}
  \caption{The  MathTOUCH equation editor}\label{fig:MathTOUCH}
\end{figure}

Equation editors which output a text-representation into a text-document, as did DragMath, will remain very helpful tools.
Such tools are especially useful for complex expressions, larger matrices and where specific symbols are not readily available.
Text-based representations can then be manipulated, e.g. parts of an expression can be copy-pasted elsewhere in the document.
Ideally, the equation editor would read in a text-representations of an existing expression (where they are syntactically correct at least) for editing.
A notable example is MathTOUCH \citep{Fukui2025}, as shown in Figure \ref{fig:MathTOUCH}. 
In MathTOUCH an AI is used to provide students with options for the likely meaning of their typed linear expression.

There are some interesting alternatives.
For example, the FlickMath interface for mobile phones uses the system for typing Japanese characters \citep{Nakamura2016b}.

The distinction between a pop-up equation editor, and inline editing of equations, is not clear cut.
There are many examples of inline-editing.
For example, MathQuill is a widely used editor, which sits inside a web page.
Users edit equations inline, within a text document, and they see and edit the immediate preview.
{\tt http://mathquill.com/} (March 2026).
Internally, all such mechanisms store the mathematics as some kind of digital object, probably encoded as a text string.
LaTeX is a popular encoding, when only presentation is needed.

\subsection{Structured proofs} \label{sec:Back}

So far we have considered the issue of typing mathematical expressions, e.g.~individual algebraic formulae.  
Now we need to consider how the semantic structure of a proof can be written effectively as well.
Written proofs serve different goals:  some proofs communicate ideas more effectively, whereas some provide rigorous justification.
Structured proofs with narrative; summary/outline; and collapsible sections provide a bridge between formally checkable proofs of Lean and informal proof outlines.
\begin{quote}
Avoiding mistakes when manipulating formulas requires careful, detailed calculations.
Avoiding mistakes when proving theorems requires careful, detailed proofs.
When first shown a detailed, structured proof, most mathematicians react: ``I don't
want to read all those details; I want to read only the general outline and perhaps
some of the more interesting parts.'' My response is that this is precisely why they
want to read a hierarchically structured proof. \citep[p.601]{Lamport1995}
\end{quote}
Many authors agree that when we have a string of equalities, each should have a brief justification, as in Figure \ref{fig:Schwamb1904-pg67}.

\begin{figure}
  \begin{center}
  \includegraphics[width=0.6\textwidth]{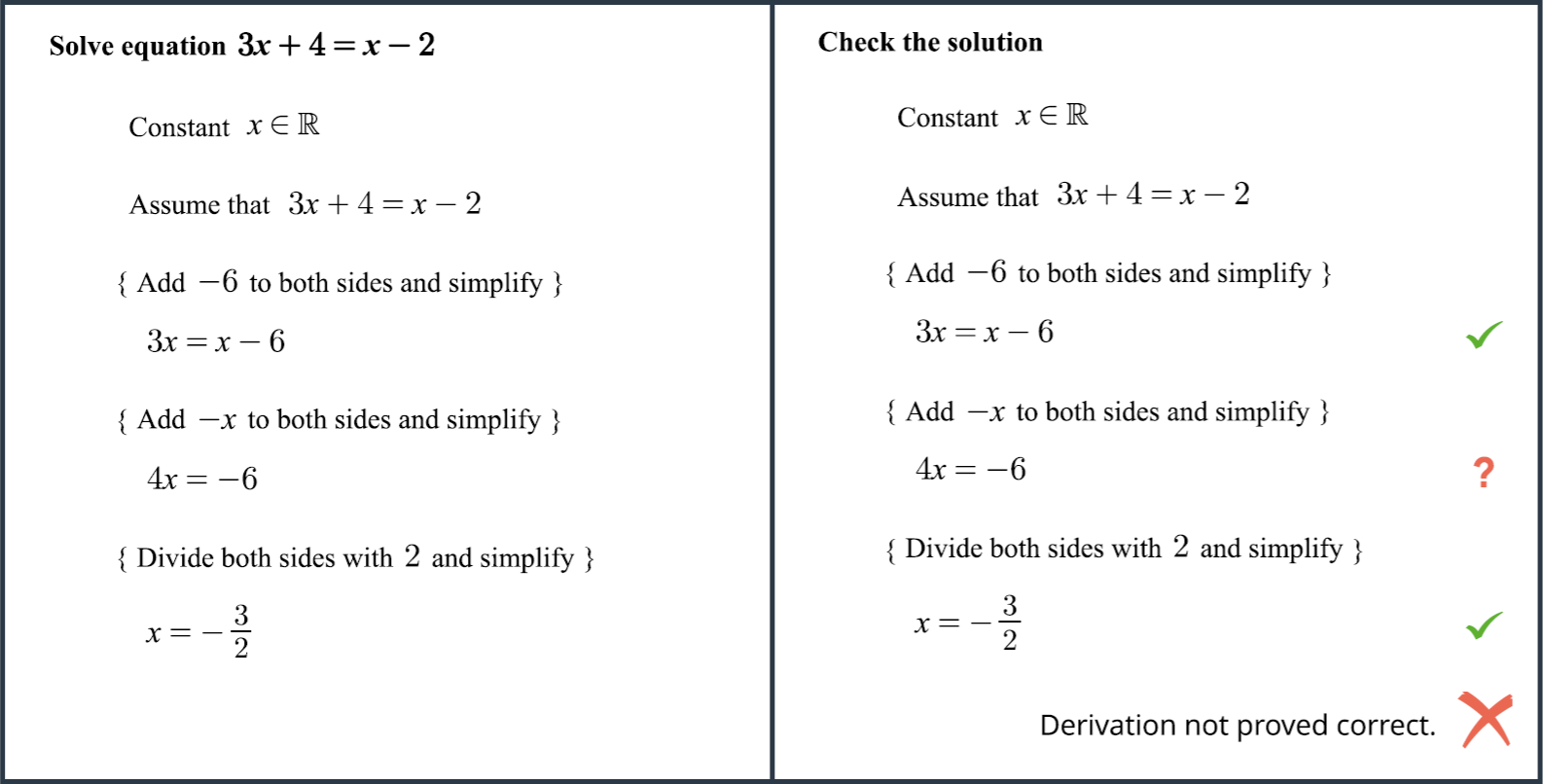}
  \end{center}
  \caption{structured derivations with eMathStudio}\label{fig:eMathStudio}
\end{figure}

There are many examples of software which aim to provide structure to mathematical arguments.
A notable example are the structured derivations suggested in \citet{Back2016}.
A basic structured derivation consists of a sequence of declaration steps, followed by a sequence of assumption steps, followed by
a sequence of conclusion steps.
An example is shown in Figure \ref{fig:eMathStudio}.
Here, the declaration is that $x$ is real, and the assumption step says that $x$ satisfies the given equation.
Three conclusion steps, with justification, complete the structured derivation.
The eMathStudio also contains some automatic assessment, with feedback, as shown on the right of Figure \ref{fig:eMathStudio}.
Further, eMathStudio has support for nested sub-proofs, which is an important feature of mathematical proof.
With these systems, there is always a compromise between freedom of expression (as provided by free text) and constrained semantic markup.
See {\tt https://emathstudio.com/} (Feb 2026).

\section{Writing in Online Mathematics Examinations in 2026}


Our goal is for students to type a short, complete mathematical argument as free-form typed text as they derive it, with minimum friction from the computer interface.
Options for diagrams, graphs and tables of information is an entirely separate issue, and are not discussed here.
Similarly, we expect a tablet for writing with OCR could contribute {\em text} as an input device, just as a keyboard does.

The system should be totally forgiving of mathematical mistakes/errors.
The system should be reasonably forgiving of errors in general, but acknowledges that students should learn how to use a tool itself.
It's reasonable to expect students to mark up the start and end of an equation.
For example, using \verb$\(...\)$ in LaTeX or \verb$`...`$ in AsciiMath.
If students forget to mark up equations, mismatch brackets, or make other mistakes then the text can still be read without (complete) loss of meaning.

The goal does not include automatic assessment of any content.
We {\em leave open the option} for future changes which would allow sections of work to be automatically checked.
For example, final algebraic answers and reasoning by equivalence of algebraic derivations can be automatically assessed, but in most work these calculations are embedded in a narrative containing justifications.
Ultimately, if the community develop proof checkers aimed at school level mathematics, these can be applied to mathematics using the proposed format.
A typed input system which includes more semantics only increases the reliability of any automatic assessment.

Our comprehensive survey of current options suggests that in 2026 it is completely practical to implement a text-based system, with immediate preview, based around the following core system:
\begin{enumerate}
  \item {\em Students type unicode text}\\
  This works cross platform, on all hardware.  It's less than ideal, but still possible, on mobile phones.
  Text-formats give complete transparency over what the student types and are future-proof.
  Proprietary formats, particularly embedded equations, are fraught with problems including obsolescence.
  Simple expressions can be native unicode.
  Text-formats enable existing accessibility tools to interact with minimum friction.
  E.g.~students can re-size the text and change text/background colour themselves.
  Existing screen readers can use the text, and speech to text systems can write plain text easily.
  Text gives transparency over what was actually typed by students.
  \item {\em Students use Markdown for document structure and formatting}\\
  Markdown is a simple, unobtrusive system for markup of simple document structures, including headings, lists, and emphasis/bold. 
  Markdown is widely used, and popular and there is excellent support for immediate markdown preview.
  \item {\em Students can include Space Math for mathematics}\\
  Space Math is ideal for school-level mathematics, and of the text-formats available is the simplest.
  It's back-compatible with most of AsciiMath, but has the added advantage of including more semantics which is very helpful.
  Note, computer code is denoted in markdown using backticks: \verb$`....`$ conflicting with the default for the default for AsciiMath.
  Space Math solves this problem by overloading the LaTeX delimiters and deciding if a particular expression is LaTeX or Space Math.
  E.g.~if a mathematical expression contains a backslash $\backslash$ or uses curly brackets with powers, e.g. \verb!x^{23}!, then the expression is treated as LaTeX.
  Otherwise, it is interpreted as Space Math.  Expressions such as \verb!$x^n$! is identical in both LaTeX and Space Math, and \verb!$x^23$! is considered as Space Math ($x^{23}$), rather than LaTeX ($x^23$).  
  \item {\em Students optionally include LaTeX mathematics environments}\\
  For short answers, there is no need for the document structure commands provided by LaTeX, but the mathematics environment is comprehensive and reliably stable.
  Many systems already embed just the LaTeX mathematics environment, and this is useful. 
  Therefore, LaTeX mathematics is optional and sits alongside Space Math without friction.
  \item {\em Equation editors}\\
  Graphical equation editors can work alongside this system.
  For example, a graphical equation editor could write Space Math into the student's document directly as plain text.
  More work would be needed to parse a student's mathematics into the editor, but that is possible with semantic markup.
\end{enumerate}

One option is to implement the above features as a twin-pane system.
Students type text into one pane (e.g.~the left) and see an immediate preview of their text in the other pane.
Adopting a text-format allows students to cut and paste between expressions.
The immediate preview shows the student all mathematics using the traditional layout of symbols, two-dimensional structures, etc.
Cut and paste provides students with the option to re-draft and edit their text, without messy crossing out of rough working.
This may prove more efficient than traditional writing, and does not interrupt the flow of thought in the same way changing to a mouse-driven equation editor seems to do.

\begin{figure}
  \begin{center}
  \includegraphics[width=0.95\textwidth]{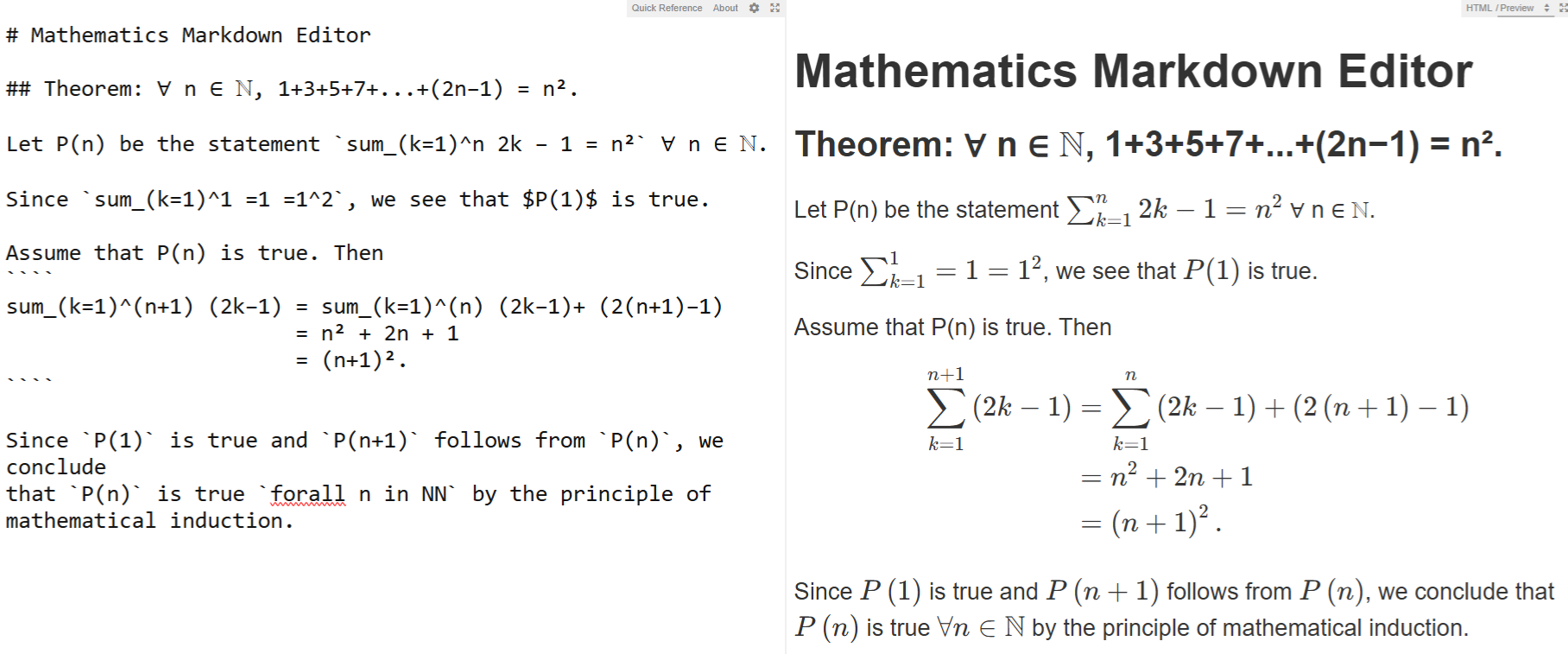}
  \end{center}
  \caption{The Minimalist Online Markdown Editor}\label{fig:2025-minimal-markdown}
\end{figure}

An example is shown in Figure \ref{fig:2025-minimal-markdown}.
On the left is the same proof shown in Figure \ref{fig:2025-unicode}, with some added markdown commands, e.g.~\verb$#$ at the start of a line to denote a heading.
Both AsciiMath and LaTeX mathematics environments have been incorporated as can be seen in the Figure.
More work is needed to include Space Math, an equation editor, or other tools.
Figure \ref{fig:2025-minimal-markdown} is to demonstrate the anticipated interface, as a practical proof of concept and for further testing and evaluation.
While, for examinations, we need support for short focused fragments of text, markdown with LaTeX has been already expanded to provide support for larger document structures in Typst {\tt https://github.com/typst/typst} (Feb 2026)
Examples such as that in Figure \ref{fig:2025-minimal-markdown} act as inspiration for rapid development and evaluation cycles for practical input systems, e.g.~in the STACK online assessment system.

Alternatively, we could implement the above features as a WYSIWYG system with input of mathematics typed by an interface similar in design to mathquill, so that equations are displayed immediately and inline.
This would parse the Space Math, as suggested, and include LaTeX equations and markdown.
Note, underneath any WYSIWYG interface would still be the markdown extended as suggested.
This gives clarity about what a student really typed and avoids propitiatory binary formats.

Students, typically, take years to learn algebra at school.  It seems reasonable to expect them to invest a little time in learning how to type, and how to type mathematics.  
If a text-based system genuinely becomes a useful tool, then this investment will be repaid.
That said, it's equally unreasonable to expect a useful tool to be {\em immediately} beneficial without some period of training and practice.

The system should contain {\em extensible blocks} which add optional functionality. Indeed LaTeX/Space Math already extends the core markdown format.
Other extensions include the following.
\begin{enumerate}
  \item {\em Extended markdown for proof structures}\\
  Markdown can be extended with specific environments for proof structures.  Figure \ref{fig:2025-minimal-markdown} shows a line by line derivation, where the lines are aligned on the $=$ signs.
  Extended markdown environments could support a variety of such structures, both presentational or semantic.
  There is a complexity cost/benefit compromise to decide how many are useful.
  PreTeX markdown already proposes such environments. 
  \item {\em Calculator tool}\\
  It seems incongruous for students to reach for a pocket calculator, or type digits into a calculator tool, and then copy the answer back into their digital document. 
  Instead, the proposal is that students can embed a calculation within their document.
  For example, the STACK system implements ``CAS-text'' which is a text-based format into which computer algebra system (CAS) calculations can be embedded.
  CAS calculations are replaced in the immediate preview by their value.
  In STACK, the syntax for such calculations is \verb${@expr@}$ where \verb$expr$ is evaluated by a CAS (e.g.~Maxima) and the tag replaced by the LaTeX expression representing the value.  For example, the text 
  \begin{verbatim}We calculate \(6\times 3={@6*3@}\).\end{verbatim}
  is displayed in the immediate preview as ``Calculate \(6\times 3 = 18\).''
  This tool is optional, enabling calculator-free examinations as required.
  If available, the tool can also be limited to numerical, scientific or computer algebra calculations.
  \item {\em Graph plotting/interactive geometry tools}\\
  Most of the diagrams in current school examinations involve $x-y$ graphs of real functions.
  Just as calculations can be embedded as CAS-text, so can traditional plots.
  For example, in STACK the command \verb${@plot(x^2/(1+x^2),[x,-3,3])@}$ is replaced by an SVG illustrating the plot of $\frac{x^2}{1+x^2}$ between $-3$ and $3$.
  Clearly, to use this kind of syntax-driven plot students need to learn the syntax.
  Just as with equation editors, graphical plot tools could be embedded.
  \item {\em Automatic assessment tool}\\
  The block-based design allows for future embedding of automatic assessment.
  The final answer could be immediately assessed (with tools available today) by having simple extensions to markdown to indicate the students' final answer.  
  For example, if a line starts with \verb$answer:expr$ then \verb$expr$ would be interpreted by existing automatic assessment systems.
  Input for multi-part questions is more complex, but students need not start with a blank document.
  A document seeded with input blocks around which students type is perfectly feasible today.
  Similarly, we can envisage an interface with other tools, such as a proof assistant or for other disciplines such as code-checkers for computing examinations (see \citet{Lobb2016}).
  Eventually we expect automatic assessment using large language models.
\end{enumerate}

One significant advantage of extending markdown in this way is the flexibility this provides for the needs of other disciplines and for future unknown needs.
For relatively short examination answers the computation burden is typically very low, and can be readily done client-side (e.g.~with Javascript) on the student's machine without significant latency.

\section{Technology and the future of mathematical activity}\label{sec:future}

The focus of this paper is to consider how simple, and traditional, elementary mathematical arguments can be typed into a computer, especially for potential use in examinations.
Supporting typed mathematics, with immediate preview and interpretation of semantic meaning, requires new technology.
In parallel we are, in 2026, seeing rapid development in large language models (LLMs) which must also be considered.

Technology has regularly affected mathematics and mathematics education for more than four hundred years.
A slide rule consists of two adjacent logarithmic scales and slide rules were routinely used by engineers for about three hundred and fifty years until the widespread introduction of cheap electronic calculators.
Invention of the slide rule involved a controversy between two rivals:  William Oughtred (1573--1660) (see \citet{Cajori1916}) and Richard Delamain (1610--1645) (see \citet[\S 122]{Taylor}).
They argued about the extent you need to understand basic principles to make effective use of a tool, a discussion still with us today.
This was widely discussed when electroinc calculators were introduced, with most educators concluding that extreme positions represent a false dichotomy.
\begin{quote}
\S378 We wish to stress that the availability of a calculator in no way reduces the need for mathematical understanding on the part of the person who is using it.
\citep[p.~112]{Cockroft1982}    
\end{quote}
Calculators have become universally adopted, with most educational jurisdictions opting for some assessments which expect students to make use of a scientific or graphical calculator and some calculator-free assessments.
Examinations in which students use computer algebra are comparatively rare.



\begin{figure}
\begin{center}
\includegraphics[width=9cm]{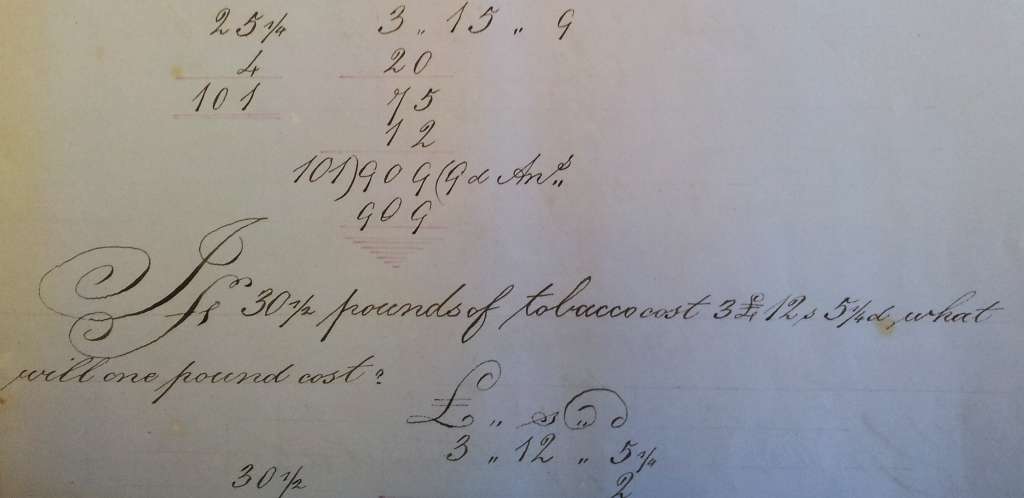}
\end{center}
\caption{Joseph Phillips's school exercise book}\label{fig:Bookwork}
\end{figure}

Let us consider Figure \ref{fig:Bookwork}, which shows a page taken from a school exercise book, attributed to Joseph Phillips, aged 10 in 1858.
Joseph is writing out long division examples, all related to commerce.
His handwriting is beautiful, and his long division is set out carefully in a standard format without any crossing out or correction.
Presumably the teacher would have assessed this work, but if so what is being assessed?
Handwriting, accountancy book-keeping practice, computational skill, or something else?
Handwriting might once have been an important part of effective communication and therefore have been assessed, however few now rely on handwriting for professional communication.
Accountancy typically now requires spreadsheets rather than paper book-keeping. 
In Section \ref{sec:tex} we mentioned how both Babbage's mechanical computers, and the introduction of TeX {\em looks backwards}, replicating existing practices and systems.
Problems arise when technology replaces something completely and renders an activity obsolete.
The computer used to be a professional occupation, but we suspect that no one is employed as a human computer today.

\citep{Hanisch1992} argued that introduction of educational technology typically follows the four phases of the {\em Didactic Inertia Principle}:
ignore; forbid; reluctantly accept, then require
\footnote{%
Prof.~Dr.Hans-Georg Weigand 
discussed this at a conference talk in 2023, and attributed the phases to \citet{Hanisch1992}.
\begin{quote}
Als Beleg für dieses Prinzip möge man sich zurückerinnem, wie Unterrichtsverwaltung und Mathematiklehrer(innen) auf den Taschenrechner reagiert haben. Dieser wurde
\begin{itemize}
  \item zuerst totgeschwiegen (Stufe 1),
  \item dann verboten (Stufe 2),
  \item mit Widerwillen erlaubt (wenn ... und aber ...) (Stufe 3) und
  \item schließlich verpfichtend eingeführt (in Österreich zur Zeit ab der 7. Schulstufe) (Stufe 4).
\end{itemize}
\citet[p.~15]{Hanisch1992}
\end{quote}
We translate this as follows:
\begin{quote}
    As evidence of this principle, consider how education administrators and math teachers reacted to the calculator. It was first ignored (Level 1), then banned (Level 2), allowed with reluctance (if... and but...) (Level 3), and finally made mandatory (in Austria, currently starting in grade 7) (Level 4).
\end{quote}
}.
Sometimes technology forces a significant change in the nature of activity.
Mathematics includes two particularly important strands of individual activity.
\begin{itemize}
  \item Developing personal skills.  For example, being able to solve problems for oneself.
  \item Learning coherent bodies of organised knowledge.  For example, calculus developed as coherent techniques rather than tricks for solving individual problems. 
\end{itemize}
These activities are related: one purpose of organising knowledge is to be able to better solve problems by finding connections between apparently different situations.

In addition to these strands of activity, students have various motivations, and these two are particularly important.
\begin{enumerate}
  \item Use in future {\em professional practice}.  
  \item {\em Cultural significance} and membership of communities of practice.  
\end{enumerate}
%
We anticipate that LLMs will soon be able to correctly solve any traditional undergraduate mathematics examination problem to a first class standard, and near perfectly in most cases.
It is speculation whether this premise is correct, or if LLMs fall slightly short of this assumption.
In 2026 LLMs can already solve many mathematics problems remarkably accurately  \citep{Ringer2025} although they continue to make basic mistakes as well.
As with other technology, it's clear that LLMs will change the nature of mathematical work, and the role humans continue to play.

\begin{figure}
  \begin{center}
  \includegraphics[width=0.85\textwidth]{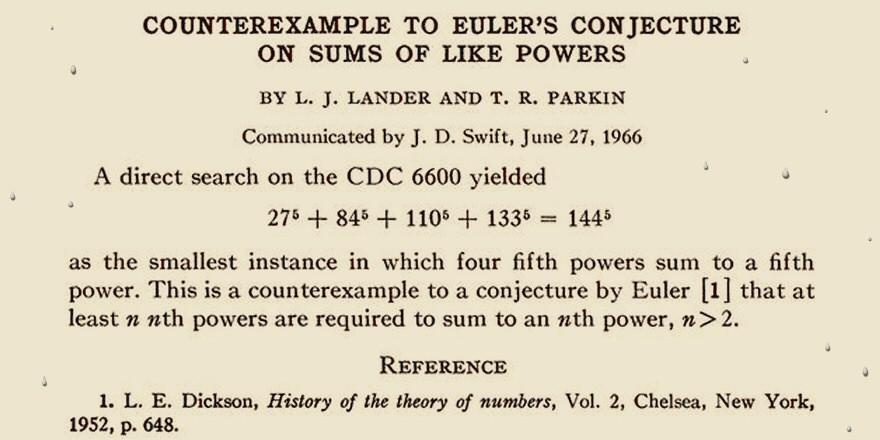}
  \end{center}
  \caption{The shortest paper in pure mathematics: \citet{Lander1967}}\label{fig:Lander1967}
\end{figure}

So, what would we like students to be able to do?
For example, how good do we expect a typical student to be at integer arithmetic?
Should students know that $25=5^2$, that $91$ is not prime, but $97$ is prime?
Obviously, there is no single answer to this question, rather each student, teacher and institution needs to answer it for themselves.
We probably don't need students to notice
\[  61917364224 = 144^5. \]
Similarly, it's almost ludicrous to suggest contemporary students could be expected to `notice' that $1729=12^3+1^3=9^3+10^3$ as Ramanjuan did \citep{Kanigel1991}.
However, despite his own prodigious calculating powers Ramanjuan acknowledges making use of the work of the professional calculator Percy Alexander MacMahon in developing his own conjectures \citep{Ramanujan1919}.
Indeed the role of the importance of calculation and exhaustive calculation in establishing the calculation
\[  27^5+84^5+110^5+133^5=144^5\]
is the subject of the shortest paper in pure mathematics \citep{Lander1967}, reproduced in Figure \ref{fig:Lander1967}.
In short, the human ability to calculate is important for the discipline of mathematics, and is likely to remain so.

Becoming part of a community is a complex process.
Humour plays a role in determining who is, and who is not, part of a community \citep{Renteln2005}.
Achievement is much more important.
Learning a classical canon of organised knowledge, including algebra and calculus, is part of demonstrating achievement.
Students might learn, by heart, proofs of some of the classical theorems just as one might choose to learn, by heart, a poem. 
E.g., one might learn one or more proofs that the square root of two is irrational \citep{Duchene2010} 
or the Pythagorean Theorem, \citep{Loomis1968}, even though a proof is little or no use when applying a theorem.

Another way students can demonstrate their achievement is solving, for themselves, problems and puzzles.
The point is that students choose to struggle to solve puzzles for themselves regardless of whether an AI could provide a flawless solution.
They gain personal satisfaction, and respect from existing members of the mathematical community, for undertaking this classical apprenticeship.
Problems have an interesting and complex history of their own.
\citet[p.~12]{Swetz2012} provides a history of problems, arguing they demonstrate the importance of classical problems in mathematics education.
Generations of students have solved, for themselves, identical problems which the mathematical community has chosen and passed down.
In this sense mathematical problems, at all levels, also play a role as cultural artefact in mathematics.

For this reason, even if LLMs {\em could} solve most problems flawlessly we expect many people will still choose to engage in the culture of mathematics.
\citet{Burkhardt2012b} called for {\em tests worth teaching to}.
Moving current examinations online, without reassessing the content and purpose, is a missed opportunity.
There are compelling examples of assessment-driven education.
R.~P.~Burn wrote three books of exercises (e.g.~\citet{Burn1987}, \citet{Burn1996} and  \citet{Burn2000}), which constitute whole university courses in which students solve problems rather than read notes.
Multiple choice questions alone were sufficient for \citet{Crowder1960} to create non-linear learning materials in undergraduate pure mathematics.
Online assessment now provides much more sophisticated tools for assessment-driven courses \citep{Kinnear2022}.
Assessments can be crafted into interesting and enjoyable experiences for students.
Indeed, practice of technical exercises need not be dull and repetitive either.
Just as musicians have {\em etudes}, designed to be a practice piece for developing a specific musical skill, so do mathematicians \citep{Foster2013}.

When navigating the disruption caused by introduction of new technology, we typically encounter four phases: ignore; forbid; reluctantly accept, then require.
To navigate these phases successfully we need to remain mindful of motivation for engaging in study, balancing any need for professional training with a need to preserve and develop culture.
Underlying current difficulties with AI is deciding what we actually want the students to be able to do without using AI (forbid) or by using AI effectively (require).
Then we need to reliably identify what our student is doing for themselves.
We need to balance personal skill development with study of pre-existing bodies of organised knowledge, perhaps with other more general skills including teamwork.
The value of education is ultimately judged by what students can do as a result, and assessments are the means by which we establish what students can do.
The tools which students have available, especially during examinations, will have a profound effect on their experiences of and dispositions towards mathematics.
The constraints of a typed format have both benefits and drawbacks.
The constraints of semantic markup enable meaning to be conveyed effectively and give the potential for automatic assessment.
Indeed, such constraints may have some benefits for beginners.
The obvious drawback is the loss of freedom to write on paper to create diagrams and annotations and to invent new notation.
Unanticipated outcomes are also likely to occur with a typed digital format.
Moving examinations on-screen effectively is therefore an important design challenge and responsibility to future generations.

\subsection*{Disclosure statement}

The second author instigated the STACK online assessment system, which is likely to incorporate ideas and results from this survey. 

\printbibliography{}

\end{document}